\documentclass[12pt]{amsart}
\usepackage{amssymb, amsmath, amsthm,enumerate,mathrsfs}

\newtheorem{theorem}{Theorem}[section]

\newtheorem{lemma}[theorem]{Lemma}

\newtheorem{prop}[theorem]{Proposition}
\theoremstyle{definition}

\newtheorem{rem}[theorem]{Remark}

\newtheorem{example}[theorem]{Example}

\def\beq{\begin{equation}
}
\def\eeq{\end{equation}}

\def\tm{\tilde{m}}

\numberwithin{equation}{section}

\newcommand{\monsub}[3]{{#1} |_{#2 \rightarrow #3}}

\def\beq{\begin{equation}}
\def\eeq{\end{equation}}

\def\cH{ {{\mathcal H}}}
\def\cN{ {{\mathcal N}}}
\def\cP{ {{\mathcal P}}}
\def\cQ{ {{\mathcal Q}}}

\def\bbR{ {\mathbb R}}

\def\gtupn{(\mathbb{R}^{n\times n})^g}

\def\cB{ {\mathcal B} }

\def\cF{ {\mathcal F} }
\def\cG{ {\mathcal G} }
\def\cH{ {\mathcal H} }

\def\cK{ {\mathcal K} }

\def\cN{ {\mathcal N} }

\def\cP{{\mathcal P}}

\def\beq{\begin{equation}}
\def\eeq{\end{equation}}

\def\cG{{\mathcal G}}

\title[Noncommutative Plurisubharmonic Polynomials: Part II]{Noncommutative Plurisubharmonic Polynomials \\ Part II: Local Assumptions}

\author[Greene]{Jeremy M. Greene${}^1$${}^\dagger$}
\address{Jeremy M. Greene, Department of Mathematics \\
University of California \\
San Diego
}
\email{j1greene@math.ucsd.edu}
\thanks{${}^1$Research supported by NSF grants
DMS-0700758, DMS-0757212, and the Ford Motor Co.}
\thanks{${}^\dagger$The material in this paper is part of the 
Ph.D. thesis of Jeremy M. Greene at UCSD}


  \makeatletter
    \newcommand{\mycontentsbox}{%
        {\centerline{NOT FOR PUBLICATION}
                  \tableofcontents}}
                    \def\enddoc@text{\ifx\@empty\@translators \else\@settranslators\fi
                        \ifx\@empty\addresses \else\@setaddresses\fi
                            \newpage\mycontentsbox}
                              \makeatother

\begin{document}

\maketitle

\begin{abstract}
We say that a symmetric noncommutative (nc) polynomial is nc plurisubharmonic 
(or nc plush) on 
an nc open set if it has an nc complex hessian that is positive semidefinite when 
evaluated on open sets of matrix tuples of sufficiently large size. 
In this paper, we show that if an nc polynomial 
is nc plurisubharmonic on an nc open set then the polynomial is actually 
nc plurisubharmonic everywhere and has the form
\beq \label{eq:MAINpsh}
p = \sum f_j^T f_j + \sum k_j k_j^T + F + F^T
\eeq
where the sums are finite and $f_j$, $k_j$, $F$ are all nc analytic.

In \cite{GHVpreprint}, it is shown that if $p$ is nc plurisubharmonic everywhere then 
$p$ has the form in Equation \eqref{eq:MAINpsh}. 
In other words, \cite{GHVpreprint} makes a global assumption while the 
current paper makes a local assumption, but both reach the same conclusion. 
We show that if $p$ is nc plurisubharmonic on an nc ``open set'' (local) then $p$ is, 
in fact, nc plurisubharmonic everywhere (global) and has the form expressed in Equation \eqref{eq:MAINpsh}. 

This paper requires a technique that is not used in \cite{GHVpreprint}. 
We use a Gram-like vector and matrix representation (called the border vector 
and middle matrix) for homogeneous degree 2 nc polynomials. 
We then analyze this representation for the nc complex hessian on an nc open set and positive semidefiniteness forces a very rigid structure on the border vector and middle matrix. 
This rigid structure plus the theorems in \cite{GHVpreprint} ultimately force 
the form in Equation \eqref{eq:MAINpsh}.
\end{abstract}

\maketitle

\section{Introduction} \label{sec:Intro}

In \cite{GHVpreprint}, it is shown that if $p$ is nc plush everywhere then 
$p$ has the form in Equation \eqref{eq:MAINpsh}. In this paper, we prove 
a stronger result on a ``local implies global'' level. We show that if $p$ is 
nc plush on an nc ``open set'' (local) then $p$ is, in fact, nc plush everywhere (global)
and has the form expressed in Equation \eqref{eq:MAINpsh}. Since this paper 
is a close companion of \cite{GHVpreprint}, we refer the reader there for 
background.

This paper requires a technique that is not used in \cite{GHVpreprint}. 
We use a Gram-like vector and matrix representation (called the border vector 
and middle matrix) for homogeneous degree 2 nc polynomials. 
We then analyze this representation for the nc complex hessian on an nc open set and positive semidefiniteness forces a very rigid structure on the border vector and middle matrix. 
This rigid structure plus the theorems in \cite{GHVpreprint} ultimately force 
the form in Equation \eqref{eq:MAINpsh}. 

\subsection{NC Polynomials, Their Derivatives, and Plurisubharmonicity} \label{subsec:BasicDefns}
This subsection and the next (Subsection \ref{subsec:differentiation}) are 
contained in Section 1 of \cite{GHVpreprint}.

\subsubsection{NC Variables and Monomials} \label{subsubsec:NCVarsMons}
We consider the free semi-group on the $2g$ noncommuting
formal variables $x_1, \ldots, x_g, x_1^T, \ldots, x_g^T$.
The variables $x_j^T$ are the formal transposes of the
variables $x_j$.  The free semi-group in these $2g$
variables generates monomials in all of these variables
$x_1, \ldots, x_g, x_1^T, \ldots, x_g^T$, often called
monomials in $x$, $x^T$. \\
\indent If $m$ is a monomial, then
$m^T$ denotes the transpose of the monomial $m$.
For example, given the monomial (in the $x_j$'s)
$x^w=x_{j_1}x_{j_2}\ldots x_{j_n}$,
the involution applied to $x^w$ is
$(x^w)^T=x_{j_n}^T \ldots x_{j_2}^T
x_{j_1}^T$.

\subsubsection{The Ring of NC Polynomials} \label{subsubsec:NCPolys}
Let $\mathbb{R} \langle x_1,\ldots,x_g, x_1^T, \ldots, x_g^T \rangle$ denote
the ring of noncommutative polynomials over $\mathbb{R}$ in the noncommuting
variables $x_1,\ldots,x_g, x_1^T,\ldots, x_g^T$.  We often abbreviate
$$
\mathbb{R} \langle x_1,\ldots,x_g,x_1^T, \ldots, x_g^T \rangle \ \ \
\text{by} \ \ \ \mathbb{R} \langle x,x^T \rangle.
$$
Note that $\mathbb{R}\langle x,x^T \rangle$ maps to itself under the
involution ${}^T$. 

We call a polynomial \textbf{nc analytic} if it contains only the
variables $x_j$ and none of the transposed variables $x_i^T$. Similarly,
we call a polynomial \textbf{nc antianalytic} if it contains only the
variables $x_j^T$ and none of the variables $x_i$.

We call an nc polynomial, $p$, \textbf{symmetric} if $p^T = p$. For
example, $p = x_1 x_1^T + x_2^T x_2$ is symmetric. The
polynomial $\tilde{p} = x_1 x_2 x_4 + x_3 x_1$ is nc analytic but not symmetric.
Finally, the polynomial $\hat{p} = x_2^T x_1^T + 4x_3^T$ is nc antianalytic
but not symmetric.

\subsubsection{Substituting Matrices for NC Variables} \label{subsubsec:SubMats}
If $p$ is an nc polynomial in the variables $x_1,\dots,x_g,x_1^T,\ldots,x_g^T$
and
$$
X=(X_1,X_2,\ldots,X_g) \in (\bbR^{n\times n})^{g},
$$
the evaluation $p(X,X^T)$ is defined by simply replacing $x_j$ by $X_j$ and $x_j^T$ by
$X_j^T$. Note that, for $Z_n=(0_n,0_n,\dots,0_n) \in (\bbR^{n\times n})^{2g}$ where
each $0_n$ is the $n\times n$ zero matrix, $p(0_n)=I_n \otimes p(0_1)$. Because of
this simple relationship, we often simply write $p(0)$ with the size $n$ unspecified.
The involution, ${}^T$, is compatible with matrix transposition, i.e.,
$$
p^T(X,X^T) = p(X,X^T)^T.
$$

\subsubsection{Matrix Positivity} \label{subsubsec:MatrixPositivity} 
We say that an nc symmetric polynomial, $p$, in the $2g$ variables 
$x_1, \ldots, x_g, x_1^T, \ldots, x_g^T$, is \textbf{matrix positive} if 
$p(X,X^T)$ is a positive semidefinite matrix when 
evaluated on every $X \in \gtupn$ for every size $n \ge 1$; i.e., 
$$
p(X,X^T) \succeq 0
$$
for all $X \in \gtupn$ and all $n \ge 1$.

\subsection{NC Differentiation} \label{subsec:differentiation}
First we make some definitions and state some properties about nc differentiation.

\subsubsection{Definition of Directional Derivative} \label{subsubsec:dird}
Let $p$ be an nc polynomial in the nc variables $x = (x_1,\ldots, x_g)$ and $x^T = (x_1^T, \ldots, x_g^T)$. 
In order to define a directional derivative, we first replace all $x_i^T$ by $y_i$. Then the \textbf{directional 
derivative of $p$ with respect to $x_j$ in the direction $h_j$} is 
\beq \label{eq:derivx}
p_{x_j}[h_j] := \frac{\partial p}{\partial x_j}(x,x^T)[h_j] = 
\frac{dp}{dt}(x_1, \ldots, x_j+th_j, \ldots, x_g, y_1,\ldots,y_g)|_{t=0} |_{y_i=x_i^T}.
\eeq
The \textbf{directional derivative of $p$ with respect to $x_j^T$ in the direction $k_j$} is 
\beq \label{eq:derivxT}
p_{x_j^T}[k_j] := \frac{\partial p}{\partial x_j^T}(x,x^T)[k_j] = 
\frac{dp}{dt}(x_1,\ldots,x_g, y_1, \ldots, y_j+tk_j, \ldots, y_g)|_{t=0} |_{y_i=x_i^T}.
\eeq
Often, we take $k_j=h_j^T$ in Equation \eqref{eq:derivxT} and we define
\begin{eqnarray*} 
p_x[h] &:=& \frac{\partial p}{\partial x}(x,x^T)[h] = \frac{dp}{dt}(x+th,y)|_{t=0}|_{y=x^T} = 
\sum_{i=1}^{g} {\frac{\partial p}{\partial x_i}(x,x^T)[h_i]} \\
p_{x^T}[h^T] &:=& \frac{\partial p}{\partial x^T}(x,x^T)[h^T] = \frac{dp}{dt}(x,y+tk)|_{t=0}|_{y=x^T, k=h^T} = 
\sum_{i=1}^{g} {\frac{\partial p}{\partial x_i^T}(x,x^T)[h_i^T]}.
\end{eqnarray*}
Then, we (abusively\footnote{For more detail, see \cite{GHVpreprint}. The idea for computing 
$p^{(\ell)}(x)[h]$ is that we first noncommutatively expand $p(x+th)$. Then, 
$p^{(\ell)}(x)[h]$ is the coefficient of $t^\ell$ multiplied by $\ell !$; i.e.,
$
p^{(\ell)}(x)[h] = (\ell !) (\mbox{coefficient of $t^\ell$}).
$ }) 
define the $\ell^{th}$ directional derivative of $p$ in the direction $h$ as 
$$
p^{(\ell)}(x)[h] := \frac{d^{\ell}p}{dt^{\ell}}(x+th,y+tk)|_{t=0}|_{y=x^T,k=h^T}
$$
so the first directional derivative of $p$ in the direction $h$ is 
\begin{eqnarray} \label{eq:dird}
p'(x)]h] &=& \frac{\partial p}{\partial x}(x,x^T)[h] + \frac{\partial p}{\partial x^T}(x,x^T)[h^T] \\
&=& p_x[h] + p_{x^T}[h^T].
\end{eqnarray}
It is important to note that the directional derivative is an nc polynomial that is 
homogeneous degree 1 in $h$, $h^T$.  If $p$ is symmetric, so is $p'$.

\begin{example} \label{ex:DirDMon}
Given a general monomial, with $c \in \bbR$,
$$ 
m = c x^{i_1}_{j_1} x^{i_2}_{j_2} \cdots x^{i_n}_{j_n} 
$$

where each $i_k$ is either 1 or $T$, we get that 
$$ 
m' =  c h^{i_1}_{j_1}x^{i_2}_{j_2}\cdots x^{i_n}_{j_n} + 
c x^{i_1}_{j_1}h^{i_2}_{j_2}x^{i_3}_{j_3}\cdots x^{i_n}_{j_n} + \cdots 
+ c x^{i_1}_{j_1}\cdots x^{i_{n-1}}_{j_{n-1}}h^{i_n}_{j_n}. 
$$
\end{example}

\subsubsection{Hessian and Complex Hessian} \label{subsubsec:Hessians}

Often, one is most interested in the hessian of a polynomial and its positivity, as this 
determines convexity. However, in this paper, we are most concerned with 
the \textit{complex hessian}, since it turns out to be related to ``nc analytic'' 
changes of variables.

We define the \textbf{nc complex hessian} \index{nc complex hessian}, $q(x,x^T)[h,h^T]$, of
an nc polynomial $p$ as the nc polynomial in the $4g$ variables $x = (x_1,\ldots, x_g)$, $x^T =
(x_1^T, \ldots, x_g^T)$, $h = (h_1, \ldots, h_g)$, and $h^T = (h_1^T, \ldots, h_g^T)$

\beq \label{def:complexhesspoly} \index{$q(x)[h]$}
q(x,x^T)[h,h^T] := \frac{\partial^2 p}{\partial s \partial t}(x+th, y+sk) |_{t,s=0} |_{y=x^T, k=h^T}.
\eeq
The nc complex 
hessian is actually a piece of the full \textbf{nc hessian} which is
\begin{eqnarray*}
p''(x)[h] &=& \frac{\partial^2 p}{\partial t^2}(x+th, y) |_{t=0} |_{y=x^T} +
\frac{\partial^2 p}{\partial t \partial s}(x+th, y+sk) |_{t,s=0} |_{y=x^T, k=h^T} \\
&+& \frac{\partial^2 p}{\partial s \partial t}(x+th, y+sk) |_{t,s=0} |_{y=x^T, k=h^T} +
\frac{\partial^2 p}{\partial s^2}(x, y+sk) |_{s=0} |_{y=x^T, k=h^T} \\
\iffalse
&=& \frac{\partial^2 p}{\partial t^2}(x+th, y) |_{t=0} |_{y=x^T} +
\frac{\partial^2 p}{\partial s^2}(x, y+sk) |_{s=0} |_{y=x^T, k=h^T} \\
&+& 2 \frac{\partial^2 p}{\partial s \partial t}(x+th, y+sk) |_{t,s=0} |_{y=x^T, k=h^T} \\
\fi
&=& 2q(x,x^T)[h,h^T] + \frac{\partial^2 p}{\partial t^2}(x+th, y) |_{t=0} |_{y=x^T} +
\frac{\partial^2 p}{\partial s^2}(x, y+sk) |_{s=0} |_{y=x^T, k=h^T}.
\end{eqnarray*}

We call a symmetric nc polynomial, $p$, \textbf{nc plurisubharmonic}
(or \textbf{nc plush}) if the nc complex hessian, $q$, of $p$ is matrix positive. 
In other words, we require that $q$ be positive
semidefinite when evaluated on all tuples of $n \times n$ matrices for every
size $n$; i.e.,
$$
q(X,X^T)[H,H^T] \succeq 0
$$
for all $X,H \in \gtupn$ for every $n \ge 1$. 

\subsection{Direct Sums} \label{subsec:DirectSums}
Our definition of the direct sum is the usual one, which for two 
matrices $X_1$ and $X_2$ is given by 
$$
X_1 \oplus X_2 := 
\left( \begin{array}{cc}
X_1 & 0 \\
0 & X_2 
\end{array} \right).
$$
Given a finite set of matrix tuples $\{X^1, \ldots, X^t \}$ with 
$$
X^j = \{ X_{j1}, X_{j2}, \ldots, X_{jg} \} \in (\mathbb{R}^{n_j \times n_j})^g
$$
for $j = 1,\ldots, t$, we define 
$$
\bigoplus_{j=1}^{t}{X^j} := \left \{ \bigoplus_{j=1}^{t}{X_{j1}}, \bigoplus_{j=1}^{t}{X_{j2}}, \ldots, 
\bigoplus_{j=1}^{t}{X_{jg}}  \right \}.
$$
For example, if $X^1=\{X_{11}, \ldots, X_{1g}\}$, $X^2=\{X_{21}, \ldots, X_{2g}\}$, 
and $X^3 = \{X_{31}, \ldots, X_{3g} \}$, we get 
$$
X^1 \oplus X^2 \oplus X^3 = \{ X_{11} \oplus X_{21} \oplus X_{31}, \ldots, 
X_{1g} \oplus X_{2g} \oplus X_{3g} \}.
$$

Now let 
$$
\cB = \bigcup_{n=1}^{\infty}{\cB_n}
$$
where $\cB_n \subseteq \gtupn$ for $n=1,2,\ldots$ be given. The graded set 
\textbf{$\cB$ respects direct sums} if for each finite set 
$$
\{X^1, \ldots, X^t \} \ \ \ \mbox{with} \ \ \ X^j \in \cB_{n_j} \ \ \ \mbox{and} \ \ \ n = \sum_{j=1}^{t}{n_j},
$$
\textit{with repetitions allowed}, $\oplus_{j=1}^{t}{X^j} \in \cB_n$.

\subsection{Noncommutative Open Set} \label{subsec:NCOpenSet}
A set $\cG \subseteq \cup_{n \ge 1}{\gtupn}$ is an \textbf{nc open set} if $\cG$ satisfies 
the following two conditions:
\begin{itemize}
\item[(i)] $\cG$ respects direct sums and 
\item[(ii)] there exists a positive integer $n_0$ such that if $n > n_0$, the set 
$\cG_n := \cG \cap \gtupn$ is an open set of matrix tuples.
\end{itemize}

We say that an nc polynomial, $p$, is \textbf{nc plush on an nc open set}, $\cG$, 
if the nc complex hessian, $q$, of $p$ satisfies
\beq \label{eq:PlushOnOpenSet}
q(X,X^T)[H,H^T] \succeq 0
\eeq
for all $X \in \cG$ and all $H \in (\mathbb{R}^{n \times n})^g$ for all $n \ge 1$. 

\subsection{Main Results} \label{subsec:MainResult}
As we will see, in Section \ref{sec:LDLandD}, the nc complex hessian, $q$, if 
matrix positive on an nc open set, can be factored as
\beq \label{eq:Star}
q = V(x,x^T)[h,h^T]^T L(x,x^T) D(x,x^T) L(x,x^T)^T V(x,x^T)[h,h^T]
\eeq
where $D(x,x^T)$ is a diagonal matrix, $L(x,x^T)$ is a lower triangular matrix with 
ones on the diagonal (we call this a \textbf{\textit{unit} lower triangular matrix}), 
and $V(x,x^T)[h,h^T]$ is a vector of monomials in $x,x^T,h,h^T$. 
The next theorem shows the surprising result that the diagonal matrix, $D(x,x^T)$, 
in Equation \eqref{eq:Star} does not depend on $x,x^T$ and that $L(x,x^T)$ has 
nc polynomial entries.

\begin{theorem} \label{thm:MainThmD}
If $p$ is an nc symmetric polynomial that is nc plurisubharmonic on an nc open set, then $q$, 
the nc complex hessian of $p$, can be written as 
$$
q = V(x,x^T)[h,h^T]^T L(x,x^T) D L(x,x^T)^T V(x,x^T)[h,h^T]
$$
where $V(x,x^T)[h,h^T]$ is a vector of monomials in $x,x^T,h,h^T$, 
$D = diag(d_1, d_2, \ldots, d_{\cN})$ is a positive semidefinite constant 
real matrix, and $L(x,x^T)$ is a unit lower triangular matrix with nc polynomial entries.
\end{theorem}

\proof The proof of this theorem requires the rest of this paper and 
culminates in Subsection \ref{subsec:Dresult}. \qed 

This gives rise to an extension of the main theorem from \cite{GHVpreprint}. 
In \cite{GHVpreprint}, it is shown that an nc polynomial which is nc plush 
\textit{everywhere} has the specific form given in Equation \eqref{eq:MainEqn} below. 
In this paper, Theorem \ref{thm:MainThm}, below, is a stronger, ``local implies global'', 
result in that an nc polynomial that is nc plush just on an nc open set is actually 
nc plush everywhere (and has the form in Equation \eqref{eq:MainEqn}).

\begin{theorem} \label{thm:MainThm}
If an nc symmetric polynomial, $p$, is nc plurisubharmonic on an nc open set, then $p$ is, in fact,
nc plurisubharmonic everywhere and has the form expressed in \cite{GHVpreprint}
\beq \label{eq:MainEqn}
p = \sum f_j^T f_j + \sum k_j k_j^T + F + F^T
\eeq
where the sums are finite and each $f_j$, $k_j$, and $F$ is nc analytic.
\end{theorem}

\proof That $D = D(x,x^T)$, in Theorem \ref{thm:MainThmD}, is a positive 
semidefinite constant real matrix immediately implies 
$$
q(X,X^T)[H,H^T] \succeq 0
$$
for all $X, H \in \cup_{n \ge 1} \gtupn$; that is, $p$ is nc plush at all $X \in \gtupn$. 
Consequently, Theorem 1.7 in \cite{GHVpreprint} gives that $p$ is of the desired form
$$
p = \sum {f_j^T f_j} + \sum {k_j k_j^T} + F + F^T
$$
where the sums are finite and $f_j$, $k_j$, $F$ are nc analytic. \qed 

Note that with an nc polynomial, $p$, as in Equation \eqref{eq:MainEqn}, the 
nc complex hessian, $q$, of $p$ is
\beq \label{eq:MainEqnCH}
q = \sum { (f_j^T)_{x^T}[h^T] (f_j)_{x}[h] } + \sum { (k_j)_{x}[h] (k_j^T)_{x^T}[h^T] } \ , 
\eeq
which is obviously matrix positive as it is a sum of squares. 
From Equation \eqref{eq:MainEqnCH}, we see that the nc complex hessian for an nc polynomial 
that is nc plush on an nc open set has even degree. 

\subsection{Guide to the Paper} \label{subsec:Guide}
In Section \ref{sec:BVMMF}, we introduce a Gram-like representation of nc quadratics. 
In Section \ref{sec:BVMMFq}, we study this Gram-like representation for the nc 
complex hessian and prove some properties for this representation. 
In Section \ref{sec:LDLandD}, we introduce the $LDL^T$ decomposition 
of the nc complex hessian and conclude that $D$ is constant. 

The author would like to thank J. William Helton and Victor Vinnikov for 
many plush discussions on noncommutative plurisubharmonic polynomials.

Now we define the Gram-like representation of an nc quadratic. It will be called 
the \textbf{middle matrix representation} (MMR).

\section{Middle Matrix Representation For A General NC Quadratic}\label{sec:BVMMF}

In this section, we turn to a special representation for nc symmetric quadratic polynomials. 
We represent nc quadratics in a factored form, $v^T M v$. This representation greatly facilitates the study 
of the positivity of nc quadratics by letting us study the positivity of $M$. Now we give details.

Any noncommutative symmetric polynomial, $f(x, x^T, h, h^T)$, in the variables $x = (x_1, \ldots, x_g)$,
$x^T =(x_1^T, \ldots, x_g^T)$, $h = (h_1,\ldots, h_g)$, and $h^T = (h_1^T, \ldots, h_g^T)$ that is degree $s$ in $x, x^T$ and homogeneous of degree two in $h, h^T$ admits a representation of the form
\beq \label{eq:genbvvmf}
f(x,x^T,h,h^T) = V(x,x^T)[h,h^T]^T M(x,x^T) V(x,x^T)[h,h^T]
\eeq
where $M(x,x^T)$, called the \textbf{middle matrix}, is a symmetric matrix of nc polynomials in $x, x^T$ and 
$V(x,x^T)[h,h^T]$, called the \textbf{border vector}, is given by
\beq \label{eq:genbv}
V(x,x^T)[h,h^T] =
\left( \begin{array}{c}
V_s(x,x^T)[h] \\
\vdots \\
V_0(x,x^T)[h] \\
V_s(x,x^T)[h^T] \\
\vdots \\
V_0(x,x^T)[h^T]
\end{array} \right).
\eeq
The $V_k(x,x^T)[h]$
(resp. $V_k(x,x^T)[h^T]$) are vectors of nc monomials of the form $h_j m(x,x^T)$
(resp. $h^T_j m(x,x^T)$) where $m(x,x^T)$ runs through the set of $(2g)^k$ monomials in $x,x^T$
of length $k$ for $j=1,\ldots, g$. Note that the degree of the monomials in $V_k$ is $k+1$.

We note that the vector of monomials, $V(x,x^T)[h,h^T]$, might contain monomials that 
are not required in the representation of the nc quadratic, $f$. Therefore, we can omit 
all monomials from the border vector that are not required. This gives us a \textit{minimal 
length} border vector and prevents extraneous zeros from occurring in the middle matrix. 
The next lemma, Lemma \ref{lem:uniquebvmmf}, says that a minimal length border vector 
contains distinct monomials.

\begin{lemma} \label{lem:uniquebvmmf}
If $f(x,x^T,h,h^T)$ is an nc symmetric polynomial that has a middle matrix representation, 
then there is a middle matrix representation for $f$ such that the border vector contains
distinct monomials. Here, distinct precludes one monomial being a scalar multiple of another.
\end{lemma}

\proof Suppose we have $f$ with the representation
$$
f(x,x^T,h,h^T) =
\left( \begin{array}{c}
m \\
\alpha m \\
n
\end{array} \right)^T
\left( \begin{array}{ccc}
p_{11} & p_{12} & p_{13} \\
p_{21} & p_{22} & p_{23} \\
p_{31} & p_{32} & p_{33}
\end{array} \right)
\left( \begin{array}{c}
m \\
\alpha m \\
n
\end{array} \right)
$$
with $\alpha$ a real number and $m$ and $n$ distinct monomials. Write $f$ as
\begin{eqnarray*}
f &=& m^T(p_{11} + \alpha^2 p_{22} + \alpha p_{21} + \alpha p_{12})m \\
&+& m^T(p_{13} + \alpha p_{23})n + n^T(p_{31} + \alpha p_{32})m + n^T p_{33} n
\end{eqnarray*}
which leads to the representation
$$
f(x,x^T,h,h^T) =
\left( \begin{array}{c}
m \\
n
\end{array} \right)
\left( \begin{array}{cc}
p_{11} + \alpha^2 p_{22} + \alpha p_{21} + \alpha p_{12} & p_{13} + \alpha p_{23} \\
p_{31} + \alpha p_{32} & p_{33}
\end{array} \right)
\left( \begin{array}{c}
m \\
n
\end{array} \right)
$$
which has distinct monomials in the border vector. \qed 

To aid us in the following sections, we cite a theorem (Theorem 8.3 in \cite{CHSY03} 
and Theorem 6.1 in \cite{HM04}). Note that in \cite{CHSY03}, the following theorem 
is stated for a positivity domain but the proof only uses the fact that positivity domains 
are nc open sets (satisfy the two conditions in Subsection \ref{subsec:NCOpenSet}). 
Hence, we slightly generalize the statement of the theorem to work on a more general 
nc open set as defined in Subsection \ref{subsec:NCOpenSet}.

\begin{theorem} \label{thm:positivity}
Consider a noncommutative polynomial $\cQ(x,x^T)[h,h^T]$ which is quadratic in the variables
$h, h^T$ that is defined on $\cG \subseteq \cup_{n \ge 1}{\gtupn}$. Write $\cQ(x,x^T)[h,h^T]$
in the form $\cQ(x,x^T)[h,h^T]=V(x,x^T)[h,h^T]^T M(x,x^T) V(x,x^T)[h,h^T]$. 
Suppose that the following two conditions hold:
\begin{itemize}
\item[(i)] the set $\cG$ is an nc open set as defined in Subsection \ref{subsec:NCOpenSet};
\item[(ii)] the border vector $V(x,x^T)[h,h^T]$ of the quadratic function $\cQ(x,x^T)[h,h^T]$ 
has distinct monomials.
\end{itemize}
Then, the following statements are equivalent:
\begin{itemize}
\item[(a)] $\cQ(X,X^T)[H,H^T]$ is a positive semidefinite matrix for each pair of tuples of matrices
$X$ and $H$ for which $X \in \cG$;
\item[(b)] $M(X,X^T) \succeq 0$ for all $X \in \cG$.
\end{itemize}
\end{theorem}

We will also need the following well known lemma (c.f. \cite{HM04}). 
Just for notational purposes of stating the lemma, let $\cB(\cH)^g$ 
denote all $g$-tuples of operators on $\cH$, where $\cH$ is a Hilbert space.

\begin{lemma} \label{lem:FTA}
Given $d$, there exists a Hilbert space $\cK$ of dimension $\sum_{0}^{2d}{(2g)^j}$ such 
that if $G$ is an open subset of $\cB(\cK)^g$, if $p$ has degree at most $d$, and if 
$p(X) = 0$ for all $X \in G$, then $p = 0$.
\end{lemma}

Now we proceed to study this middle matrix representation for the nc complex hessian.

\section{Middle Matrix Representation For The NC Complex Hessian} \label{sec:BVMMFq}

In Section \ref{sec:BVMMF}, we introduced the middle matrix representation 
for a general nc quadratic polynomial, and this section specializes it to the nc 
complex hessian. The requirement that the nc complex hessian be positive 
on an nc open set forces rigid structure to the border vector and middle matrix.

\subsection{Border Vector for a Complex Hessian: Choosing an Order for Monomials} \label{subsec:MonOrder}
Let $p$ be an nc symmetric polynomial in $g$ variables such that the degree of its nc 
complex hessian is $d$. Then the complex hessian will be homogeneous of degree two in $h,h^T$.

For a fixed degree $k$, there are $g^k$ nc analytic monomials and $g^k$ nc antianalytic monomials
in $x,x^T$. That means there are $(2g)^k - g^k - g^k = (2g)^k - 2g^k$ `mixed' monomials of
degree $k$ (i.e., monomials that are not nc analytic nor nc antianalytic).

\subsubsection{Analytic Border Vector} \index{analytic border vector} \label{subsec:AnalyticBV}
For $0 \le k \le d-2$, let $A_k = A_k(x)[h]$ be the vector of nc analytic monomials with entries
$h_j m(x)$ where $m(x)$ runs through the set of $g^k$ nc analytic monomials of length $k$ for
$j=1,\ldots, g$. The order we impose on the monomials in this vector is lexicographic order. 
Thus, the length of $A_k = A_k(x)[h]$ is $g^{k+1}$ and the vector
\beq \label{def:analyticbv} \index{$A(x)[h]$}
A(x)[h] = col(A_{d-2}, \ldots, A_1, A_0)
\eeq
has length $g^{d-1} + \cdots + g^2 + g = g\nu$ where $\nu = g^{d-2} + \cdots + g^2 + g + 1$.

\subsubsection{Antianalytic Border Vector} \label{subsec:AntiAnalyticBV} \index{antianalytic border vector}
Let $A_k^t = A_k(x^T)[h^T]$ be the same as $A_k = A_k(x)[h]$ except replace
each $h_j$ with $h_j^T$ and replace each $x_i$ by $x_i^T$. So $A_k^t$ is the vector
of nc antianalytic monomials with entries $h_j^T m(x^T)$ where $m(x^T)$ runs through
the set of $g^k$ nc antianalytic monomials of length $k$ for $j=1,\ldots, g$ (again,
the order is lexicographic). Thus, the length of $A_k^t = A_k(x^T)[h^T]$ is
$g^{k+1}$ and the vector
\beq \label{def:antianalyticbv} \index{$A(x^T)[h^T]$}
A(x^T)[h^T] = col(A_{d-2}^t, \ldots, A_1^t, A_0^t)
\eeq
also has length $g\nu$.

\subsubsection{Mixed Term Border Vector} \label{subsec:MixedBV} \index{mixed term border vector}
Next, we define notation to handle all nonanalytic and nonantianalytic monomials.
Let $B_1 = B_1(x,x^T)[h]$ be the vector of monomials with entries $h_j x_i^T$ for
$i = 1,\ldots, g$ and $j = 1,\ldots, g$. The length of $B_1$ is $g^2$.
For $2 \le k \le d-2$, let $B_k = B_k(x,x^T)[h]$
be the vector of monomials with entries $h_j m(x,x^T)$ where $m(x,x^T)$ runs through
the set of $(2g)^k - 2g^k$ monomials of length $k$ that are not nc analytic nor nc antianalytic
for $j=1,\ldots, g$. Again, we put the same lexicographic order on the monomials.
Thus, the length of $B_k = B_k(x,x^T)[h]$ is $g((2g)^k - 2g^k)$ and the vector
\begin{equation*} \index{$B(x,x^T)[h]$}
B(x,x^T)[h] = col(B_{d-2}, \ldots, B_2, B_1)
\end{equation*}
has length $g^2 + \sum_{k=2}^{d-2} {g((2g)^k - 2g^k)}$. Then we can also define
$B_1^t = B_1^t(x,x^T)[h^T]$ to be the vector of monomials with entries $h_j^T x_i$
for $i=1,\ldots, g$ and $j=1,\ldots, g$. This also has length $g^2$. Then we define,
for $2 \le k \le d-2$, the vector $B_k^t = B_k(x,x^T)[h^T]$ to be the same as
$B_k$ except $h_j$ is replaced by $h_j^T$.  In other words, each entry looks like
$h_j^T m(x,x^T)$. Then the vector
\begin{equation*} \index{$B(x,x^T)[h^T]$}
B(x,x^T)[h^T] = col(B_{d-2}^t, \ldots, B_2^t, B_1^t)
\end{equation*}
has the same length as $B(x,x^T)[h]$. 

Note that the degree of the monomials in $A_k, A_k^t, B_k, B_k^t$ is $k+1$.

\subsection{The Middle Matrix of a Complex Hessian} \label{subsec:MM}
Now we can represent the nc complex hessian, $q$,
of a symmetric nc polynomial $p$ as
\begin{equation} \label{eq:generalbvmmfforq}
q(x,x^T)[h,h^T]=
\left(
\begin{array}{c}
A(x)[h] \\
B(x,x^T)[h] \\
A(x^T)[h^T] \\
B(x,x^T)[h^T]
\end{array}
\right)^T
\left(
\begin{array}{cccc}
Q_1 & Q_2 & 0 & 0 \\
Q_2^T & Q_4 & 0 & 0 \\
0 & 0 & Q_5 & Q_6 \\
0 & 0 & Q_6^T & Q_8
\end{array}
\right)
\left(
\begin{array}{c}
A(x)[h] \\
B(x,x^T)[h] \\
A(x^T)[h^T] \\
B(x,x^T)[h^T]
\end{array}
\right)
\end{equation}
where $Q_i = Q_i(x,x^T)$ are matrices with nc polynomial entries in the variables
$x_1, \ldots, x_g$, $x_1^T, \ldots, x_g^T$. 

Again, we wish to stress that the vectors 
$A(x)[h]$, $A(x^T)[h^T]$, $B(x,x^T)[h]$, and $B(x,x^T)[h^T]$ 
may contain monomials that are not required in the representation of the complex 
hessian, $q$. Therefore, we omit all monomials from the border vector that 
are not required. This gives us a \textit{minimal length} border vector and prevents 
extraneous zeros from occurring in the middle matrix. Lemma \ref{lem:uniquebvmmf} 
says that a minimal length border vector contains only distinct monomials. 

The next subsection provides some necessary background on nc differentiation.

\subsection{Levi-differentially Wed Monomials} \label{subsec:Levi}
An extremely important fact about the nc complex hessian, $q(x,x^T)[h,h^T]$, 
is that it is quadratic in $h,h^T$ and that \textbf{each term contains some $h_j$ 
and some $h_k^T$}. If a certain monomial $m$ is in $q$, then any monomial obtained by exchanging 
$h_j$ with some $x_\ell$ in $m$ and/or exchanging some $h_k^T$ with some $x_j^T$ in $m$ is 
also in $q$. We say two such monomials are \textbf{Levi-differentially wed}. 
Indeed, being Levi-differentially wed is an equivalence relation on the monomials in 
$q$ with the coefficients of all Levi-differentially wed monomials in $q$ being the same.

\begin{example} \label{ex:LeviWed1}
The monomials $h^T h x^T x$, $h^T x x^T h$, $x^T h h^T x$, and $x^T x h^T h$ 
are all Levi-differentially wed to each other. 
\end{example}

\begin{example} \label{ex:LeviWed2}
None of the monomials $h^T h x^T x$, $h^T x h^T x$, $x^T h x^T h$ are 
Levi-differentially wed to each other.
\end{example}

The next theorem gives necessary and sufficient conditions as to when an nc 
polynomial is an nc complex hessian. This theorem is proved in \cite{GHVpreprint} 
but gets used extensively in this paper.

\begin{theorem} \label{thm:P1andP2}
An nc polynomial $q$ in $x,x^T,h,h^T$ is an nc complex hessian if and only 
if the following two conditions hold:

\begin{itemize}
\item[(P1)] Each monomial in $q$ contains exactly one $h_j$ 
and one $h_k^T$ for some $j,k$. 

\item[(P2)] If a certain monomial $m$ is contained in $q$, any 
monomial $\tm$ that is Levi-differentially wed to $m$ is also 
contained in $q$.
\end{itemize}
\end{theorem}

\proof The proof is provided in \cite{GHVpreprint}. \qed 

Theorem \ref{thm:P1andP2} (P1) shows that every term in the complex hessian, $q$, 
has an $h_j$ and $h_k^T$ for some $j$ and $k$. This structure forces the zeros in the 
middle matrix in Equation \eqref{eq:generalbvmmfforq} above. 

\subsection{Structure of the Middle Matrix} \label{subsec:BVMMFProperties}
In this subsection, we prove some properties about the structure of the middle 
matrix in the MMR for a matrix positive nc complex hessian.

\begin{lemma}\label{lem:factorization}
Let $p$ be an nc symmetric polynomial that is nc plush on an nc open set, $\cG$.
Then, the MMR in Equation \eqref{eq:generalbvmmfforq} for its nc complex hessian, 
$q$, of $p$ has $Q_2 = Q_4 = Q_6 = Q_8 = 0$. Thus, 
\beq \label{eq:bvmmfforq}
q =
\left( \begin{array}{c}
A(x)[h] \\
A(x^T)[h^T]
\end{array} \right)^T
\left( \begin{array}{cc}
Q_1(x,x^T) & 0  \\
0 & Q_5(x,x^T)
\end{array}\right)
\left(\begin{array}{c}
A(x)[h] \\
A(x^T)[h^T]
\end{array}\right).
\eeq
\end{lemma}

\proof
We consider the upper left block of the middle matrix in Equation \eqref{eq:generalbvmmfforq}
\begin{equation*}
\left(
\begin{array}{c}
A(x)[h] \\
B(x,x^T)[h]
\end{array}
\right)^T
\left(
\begin{array}{cc}
Q_1(x,x^T) & Q_2(x,x^T) \\
Q_2(x,x^T)^T & Q_4(x,x^T)
\end{array}
\right)
\left(
\begin{array}{c}
A(x)[h]\\
B(x,x^T)[h]
\end{array}
\right)
\end{equation*}
with the goal of showing $Q_2 = 0$ and $Q_4 = 0$.
Thus, suppose the border vector contains a nonzero 
monomial which is an entry in the vector of mixed monomials,
$B(x,x^T)[h]$; i.e., the border vector contains a term
\beq \label{eq:bv1}
h_k m_1(x,x^T) x_j^T m_2(x,x^T)
\eeq
for some monomials $m_1$ and $m_2$ in the variables $x_1, \ldots, x_g, x_1^T, \ldots, x_g^T$. 

Soon we shall look at the diagonal entry, $\cP^{(0)}$, in the middle matrix corresponding to this border
vector monomial in \eqref{eq:bv1} and show it is 0. By Theorem \ref{thm:positivity}, we have the middle matrix positive semidefinite for every $X$ in the nc open set, $\cG$. By Lemma \ref{lem:FTA}, if an 
nc polynomial is zero on an open set of matrix tuples with sufficiently large dimension, then the nc 
polynomial is identically zero. Hence, if there is ever a diagonal entry in the middle matrix that is 
zero on an open set of matrix tuples of large enough dimension, then that diagonal entry is identically 
zero. Hence, to force matrix positivity, the corresponding row and column in the middle matrix must 
be zero. This implies that the particular monomial in the border vector is not needed in the 
representation, thereby contradicting the border vector being of minimal length. Thus, showing $\cP^{(0)}$ 
is 0, a contradiction. 

The term(s) in the nc complex hessian corresponding to the diagonal entry $\cP^{(0)}$ 
of the middle matrix and monomial \eqref{eq:bv1} in the border vector are
\begin{equation*}
m_2^T x_j m_1^T h_k^T \mathcal{P}^{(0)} h_k m_1 x_j^T m_2
\end{equation*}
where $\mathcal{P}^{(0)}$ is some matrix positive polynomial in $x_1,\ldots, x_g, x_1^T, \ldots, x_g^T$. 
By Theorem \ref{thm:P1andP2} (P2), $q$ must also contain the Levi-differentially wed term(s)
\begin{equation*}
m_2^T h_j m_1^T h_k^T \mathcal{P}^{(0)} x_k m_1 x_j^T m_2.
\end{equation*}
This means the border vector must contain the monomial(s)
\beq \label{eq:bv2}
\{ h_k^T \mathcal{P}^{(0)} x_k m_1 x_j^T m_2 \}_{mon}
\eeq
where $\{ h_k^T \mathcal{P}^{(0)} x_k m_1 x_j^T m_2 \}_{mon}$ is the list of the monomials that 
appear as terms in the nc polynomial 
$h_k^T \mathcal{P}^{(0)} x_k m_1 x_j^T m_2$.

Again, we shall look at the term(s) in $q$ corresponding to the diagonal in the middle matrix 
corresponding to any one of the border vector monomial(s) in \eqref{eq:bv2}. Pick $h_k^T \widehat{\mathcal{P}}^{(0)} x_k m_1 x_j^T m_2$ as a specific border vector monomial in the list in \eqref{eq:bv2}. Then, the term(s) in $q$ look like
\begin{equation*}
m_2^T x_j m_1^T x_k^T (\widehat{\mathcal{P}}^{(0)})^T h_k \mathcal{P}^{(1)} h_k^T 
\widehat{\mathcal{P}}^{(0)} x_k m_1 x_j^T m_2
\end{equation*}
where $\mathcal{P}^{(1)}$ is a matrix positive polynomial in $x_1,\ldots, x_g, x_1^T, \ldots, x_g^T$, 
which is a diagonal entry of the middle matrix.
Theorem \ref{thm:P1andP2} (P2) implies $q$ must also contain the Levi-differentially wed term(s)
\begin{equation*}
m_2^T h_j m_1^T h_k^T (\widehat{\mathcal{P}}^{(0)})^T x_k \mathcal{P}^{(1)} x_k^T 
\widehat{\mathcal{P}}^{(0)} x_k m_1 x_j^T m_2
\end{equation*}
which means the border vector must contain the monomial(s)
\beq \label{eq:bv3}
\{ h_k^T (\widehat{\mathcal{P}}^{(0)})^T x_k \mathcal{P}^{(1)} x_k^T 
\widehat{\mathcal{P}}^{(0)} x_k m_1 x_j^T m_2 \}_{mon}
\eeq
where $\{ h_k^T (\widehat{\mathcal{P}}^{(0)})^T x_k \mathcal{P}^{(1)} x_k^T 
\widehat{\mathcal{P}}^{(0)} x_k m_1 x_j^T m_2 \}_{mon}$ is the list of the monomials that appear as terms 
in the nc polynomial $h_k^T (\widehat{\mathcal{P}}^{(0)})^T x_k \mathcal{P}^{(1)} x_k^T 
\widehat{\mathcal{P}}^{(0)} x_k m_1 x_j^T m_2$.

Note that the border vector monomial in \eqref{eq:bv3} has degree at least 2 more than the degree of the 
border vector monomial in \eqref{eq:bv2} which has degree at least 2 more than the degree of the border vector monomial in \eqref{eq:bv1}. We can continue this process and the degree of the successive border vector monomials will keep increasing by at least 2 at each step.  At some step, the degree of the border vector monomial will exceed $d-1$. This contradicts the fact that the border vector monomials must have degree at most $d-1$.  Thus, we have shown that $Q_4 = 0$. A similar argument shows that $Q_8 = 0$. Since the middle matrix is positive semidefinite, we also get $Q_2 = 0$ and $Q_6 = 0$, by the argument in the previous paragraph. Hence, the nc complex hessian has the representation in Equation 
\eqref{eq:bvmmfforq}, as claimed by the theorem. \qed

We call an nc polynomial \textbf{hereditary} \index{hereditary polynomial} if all $x_1^T, x_2^T,
\ldots x_g^T$ variables appear to the \textit{left} of every $x_1, x_2, \ldots, x_g$ variable.
Similarly, we call an nc polynomial \textbf{antihereditary} \index{antihereditary polynomial} if
all $x_1^T, x_2^T, \ldots x_g^T$ variables appear to the \textit{right} of every $x_1, x_2,
\ldots, x_g$ variable. 

\begin{theorem}\label{thm:MMheredantihered}
The nc complex hessian, $q$, of an nc symmetric polynomial that is nc plush on an nc open set
can be written as in Equation \eqref{eq:bvmmfforq} 
\begin{equation*} \index{border vector middle matrix representation for $q$}
q(x,x^T)[h,h^T] = \left( \begin{array}{c}
A(x)[h] \\
A(x^T)[h^T]
\end{array}\right)^T
\left( \begin{array}{cc}
Q_1(x,x^T) & 0 \\
0 & Q_5(x,x^T)
\end{array}\right)
\left( \begin{array}{c}
A(x)[h] \\
A(x^T)[h^T]
\end{array}\right)
\end{equation*}
where every nc polynomial entry in $Q_1(x,x^T)$ is hereditary and every nc polynomial entry in $Q_5(x,x^T)$
is antihereditary.
\end{theorem}

\proof
Suppose, for the sake of contradiction, $Q_1$ contains an nc polynomial entry which is not 
hereditary. Without loss of generality, this nc polynomial contains a term of the form
\beq \label{eq:badMMpoly}
m_1(x^T) x_j x_k^T m_2(x,x^T)
\eeq
where $m_1$ is a monomial in $x^T$ and $m_2$ is a monomial in $x$ and $x^T$.
Since this is part of an entry in the middle matrix, this means that the nc complex hessian must
contain a term of the form
\begin{equation*}
m_3(x^T) h_\ell^T m_1(x^T) x_j x_k^T m_2(x,x^T) h_s m_4(x) 
\end{equation*}
where $m_3(x^T) h_\ell^T$ is a specific monomial entry from the vector $A(x)[h]^T$ and $h_s m_4(x)$ is a specific monomial entry from the vector $A(x)[h]$. Then, Theorem \ref{thm:P1andP2} (P2) implies that the nc complex hessian must also contain the Levi-differentially wed term
\begin{equation*}
m_3(x^T) h_\ell^T m_1(x^T) h_j x_k^T m_2(x,x^T) x_s m_4(x).
\end{equation*}
This implies that the border vector must contain the monomial
\begin{equation*}
h_j x_k^T m_2(x,x^T) x_s m_4(x)
\end{equation*}
which contradicts having an nc analytic or nc antianalytic border vector, as required by 
Lemma \ref{lem:factorization}. The proof that $Q_5$ contains antihereditary nc polynomial 
entries is similar.
\qed  

For a real number, $r$, we define $\lfloor r \rfloor$ as the largest integer less than or equal to $r$ and 
we define $\lceil r \rceil$ as the smallest integer greater than or equal to $r$. The next theorem puts an upper bound on the degree of the monomials in the border vector for $q$. 

\begin{lemma} \label{lem:degbdforbv}
Suppose $p$ is an nc symmetric polynomial that is nc plush on an nc open set. If the degree
of its nc complex hessian, $q$, is $d$, then the degree of the border vector monomials is at 
most $\left\lfloor  \frac{d}{2} \right\rfloor$.
\end{lemma}

\proof Write the MMR for $q(x,x^T)[h,h^T]$ as
$$
q = V^T M V =
\left( \begin{array}{cc}
V_1^T & V_2^T
\end{array} \right)
\left( \begin{array}{cc}
M_1 & M_2 \\
M_2^T & M_4
\end{array} \right)
\left( \begin{array}{c}
V_1 \\
V_2
\end{array} \right)
$$
with the following property. If $d$ is odd, $V_1$ contains monomials of degree
$1,\ldots,\left\lfloor \frac{d}{2} \right\rfloor$ and
$V_2$ contains monomials of degree $\left\lceil \frac{d}{2} \right\rceil, \ldots, d-1$.
If $d$ is even, $V_1$ contains monomials of degree $1,\ldots,\frac{d}{2}$ and
$V_2$ contains monomials of degree $\frac{d}{2} + 1, \ldots, d-1$.  
In either case, polynomials in $M_4$ correspond to terms in $q$ having degree strictly greater
than $d$. Hence $M_4 = 0$. By Theorem \ref{thm:positivity}, $M(X) \succeq 0$ for all $X$ in an nc open set. 
This forces $M_2(X) = 0$ for all $X$ in an nc open set. Then, by taking $X$ to 
have large enough size, Lemma \ref{lem:FTA} implies $M_2 = 0$. \qed 

\subsubsection{Consequences of Positivity of the Complex Hessian} \label{subsubsec:ConsCH}

Now we turn from a description of the middle matrix to describing the structure of the 
nc complex hessian of an nc polynomial that is nc plush on an nc open set.

\begin{prop} \label{prop:heredplusantihered}
The nc complex hessian, $q$, of an nc symmetric polynomial that is nc plush on an nc open set
is a sum of hereditary and antihereditary polynomials.
\end{prop}
\proof
This follows immediately from Lemma \ref{lem:factorization} and Theorem \ref{thm:MMheredantihered}.
\qed 

Finally, we show that the degree of $q$ must be even when $p$ is nc plush on an nc open set. 
This fact is obvious if $p$ is assumed nc plush everywhere because then the nc complex hessian 
is a sum of squares.

\begin{theorem} \label{thm:qEvenDeg}
Suppose $p$ is an nc symmetric polynomial that is nc plush on an nc open set. Then, the degree of
its complex hessian, $q$, is even.
\end{theorem}

\proof Suppose the degree of $q$ is $2N+1$. Without loss of generality, 
Proposition \ref{prop:heredplusantihered} and Theorem \ref{thm:P1andP2}, 
requiring the presence of Levi-differentially wed monomials, 
imply that $q$ must contain a hereditary term of the form
$$
x^T_{i_1} x^T_{i_2} \cdots h^T_{i_s} h_{j_1} x_{j_2} \cdots x_{j_\ell}
$$
where $s, \ell > 0$, $s+\ell = 2N+1$, and $i_1, \ldots, i_s, j_1, \ldots, j_\ell \in \{1, \ldots, g\}$.
This means that in the middle matrix representation for $q$, the border vector must 
contain $h_{j_1} x_{j_2} \cdots x_{j_\ell}$
and $h_{i_s} x_{i_{s-1}} \cdots x_{i_1}$ which have degree $\ell$ and $s$, respectively.
But since $s+\ell = 2N+1$ and $s,\ell > 0$, one of either $s$ or $\ell$ is at least
$\left\lceil \frac{2N+1}{2} \right\rceil$. This contradicts Lemma \ref{lem:degbdforbv}.
\qed 

\section{$LDL^T$ Decomposition Has Constant $D$} \label{sec:LDLandD}
This section concerns the ``algebraic Cholesky'' factorization, $LDL^T$, of the middle matrix. 
We will show that for an nc polynomial that is nc plush on an nc open set, this $D$ is
a positive semidefinite matrix whose diagonal entries are all nonnegative real constants, and $L$ 
is unit lower triangular with entries which are nc polynomials. This is a stronger conclusion than one 
would expect because, typically, such factorizations have nc rational entries, see \cite{CHSY03, HMV06}. 
In our approach, the $LDL^T$ factorization of a symmetric matrix with
noncommutative entries will be the key tool for the determination of the matrix
positivity of an nc quadratic function.

\subsection{The $LDL^T$ Decomposition} \label{subsec:LDL}
Begin by considering the block $2 \times 2$ matrix
$$
M = \left( \begin{array}{cc}
A & B^T \\
B & C
\end{array} \right)
$$
where $A$ is a constant real symmetric invertible matrix and $B$ and $C$ 
are matrices with nc polynomial entries with $C$ symmetric. Then, $M$ has 
the following decomposition
\beq \label{eq:Step1LDLT}
M = 
\left( \begin{array}{cc}
I & 0 \\
B A^{-1} & I
\end{array} \right)
\left( \begin{array}{cc}
A & 0 \\
0 & C - BA^{-1}B^T 
\end{array} \right)
\left( \begin{array}{cc}
I & A^{-1}B^T \\
0 & I
\end{array}\right),
\eeq
where all matrices in this decomposition contain nc polynomial entries. 
If $C - BA^{-1}B^T$ contains a constant real symmetric invertible matrix somewhere 
on the diagonal, then we can apply a permutation, $\Pi$, on the left of $M$ and its 
transpose, $\Pi^T$, on the right of $M$ to move this constant real symmetric invertible matrix 
to the first (block) diagonal position of $C - BA^{-1}B^T$. We then pivot off this constant 
real symmetric invertible matrix, factor $C - BA^{-1}B^T$ as $\hat{L}\hat{D}\hat{L}^T$, and we get 
$$
\Pi M \Pi^T = 
\left( \begin{array}{cc}
I & 0 \\
BA^{-1} & \hat{L}
\end{array} \right)
\left( \begin{array}{cc}
A & 0 \\
0 & \hat{D} 
\end{array} \right)
\left( \begin{array}{cc}
I & A^{-1}B^T \\
0 & \hat{L}^T
\end{array} \right).
$$
This can be continued, provided at each step, a constant real symmetric invertible matrix  appears somewhere on the diagonal to obtain $\Pi M \Pi^T = L D L^T$ where $L$ is a unit lower triangular matrix with nc polynomial entries and $D$ is a (block) diagonal matrix with real constant blocks. 
This special situation is the one which turns out to hold in the derivation which follows.

Indeed, we shall only care about the case where $A$ is a constant real symmetric invertible 
matrix. For the case where $A$ contains nc polynomial entries and is considered to be 
``noncommutative invertible'', please see \cite{CHSY03}. In this case, we also have the 
notion of ``noncommutative rational'' functions (please see \cite{HMV06}). However, as we 
soon shall see, while nc rationals are mentioned, they never actually appear in any calculations 
in this paper.

We recall an immediate consequence of Theorem 3.3 in \cite{CHSY03}:

\begin{theorem} \label{thm:LDLTCHSY}
Suppose $M(x,x^T)$ is a symmetric $r \times r$ matrix with noncommutative
rational function entries and that $M(X,X^T) \succeq 0$ for all $X$ in some nc open set. 
Then, there exists a permutation matrix, $\Pi$, a diagonal matrix, $D(x,x^T)$, with nc rational entries, 
and a unit lower triangular matrix, $L(x,x^T)$, with nc rational entries such that 
$$
\Pi M(x,x^T) \Pi^T = L(x,x^T)D(x,x^T)L(x,x^T)^T.
$$
\end{theorem}

\begin{rem} \label{rem:perm}
In this paper, we care about the positivity of the middle matrix, $M(x,x^T)$. If $\Pi$ 
is a permutation matrix, it is clear that 
$$
\Pi M(X,X^T) \Pi^T \succeq 0 \ \ \ \ \ \Longleftrightarrow \ \ \ \ \ 
M(X,X^T) \succeq 0
$$ 
for any $X \in \mathbb{R}^{n \times n}$ and any $n \ge 1$. 
As a result, for ease of exposition, we will often, without loss of 
generality, omit the permutation matrix, $\Pi$. 

Also, there will be some instances where we will, without loss of generality, 
assume a specific order in the border vector, $V(x,x^T)[h,h^T]$. For example, 
we may assume a given monomial, say, $h m(x,x^T)$, is the first monomial 
in $V(x,x^T)[h,h^T]$. This assumption also amounts to a permutation of 
$V(x,x^T)[h,h^T]$ which, again, does not affect positivity of $M(x,x^T)$ 
so we omit it from the discussion.
\end{rem}

We now proceed to apply the $LDL^T$ factorization to the middle matrix of the nc complex hessian. 
Let $p$ be an nc symmetric polynomial and let $q$ denote the nc complex hessian of $p$.
Since $q$ is homogeneous of degree 2 in $h,h^T$, $q$ admits the MMR 
\beq \label{eq:basicBVMMF}
q = V(x,x^T)[h,h^T]^T M(x,x^T) V(x,x^T)[h,h^T].
\eeq
If $p$ is nc plush on an nc open set, then $M(x,x^T)$ is symmetric and matrix positive on an nc open set 
and we can factor $M(x,x^T)$ following the process underlying Equation \eqref{eq:Step1LDLT} 
and Theorem \ref{thm:LDLTCHSY}, thus converting Equation \eqref{eq:basicBVMMF} to
\beq \label{eq:qLDLT}
q = V(x,x^T)[h,h^T]^T L(x,x^T) D(x,x^T) L(x,x^T)^T V(x,x^T)[h,h^T]
\eeq
up to a harmless rearrangement of the border vector.

In Section \ref{subsec:Dresult}, we prove one of the main theorems of this paper, Theorem 
\ref{thm:ConstantD}, which was stated in Section \ref{subsec:MainResult} as Theorem 
\ref{thm:MainThmD}. We recall that this theorem says that $D(x,x^T)$ in Equation 
\eqref{eq:qLDLT} does not depend on $x,x^T$ and is a positive semidefinite 
constant real diagonal matrix for an nc polynomial that is nc plush on an nc open set. 
In addition, we will prove that $L(x,x^T)$ contains nc polynomials instead of nc rationals. 
Now we start the build up to Section \ref{subsec:Dresult}.

\subsection{Properties of $LDL^T$ for NC Polynomials that are NC Plush on an NC Open Set} \label{subsec:LDLplush}
In this subsection, we present properties of the $LDL^T$ factorization of the
nc complex hessian for an nc polynomial that is nc plush on an nc open set.

Recall from Section \ref{subsec:NCOpenSet} that a set $\cG \subseteq \cup_{n \ge 1} \gtupn$ 
is an nc open set if:
\begin{itemize}
\item[(i)] $\cG$ respects direct sums, and
\item[(ii)] there exists a positive integer $n_0$ such that if $n > n_0$, the set 
$\cG_n := \cG \cap \gtupn$ is an open set of matrix tuples;
\end{itemize}
and an nc symmetric polynomial, $p$, is nc plush on an nc open set, $\cG$, 
if $p$ has an nc complex hessian,
$q$, such that $ q(X,X^T)[H,H^T] $ is positive semidefinite for all $X \in \cG$ 
and for all $H \in (\bbR^{n\times n})^g$ for every $n \ge 1$.

Theorem \ref{thm:qEvenDeg} shows that the nc complex hessian has even degree; denote it $2N$. 
We will use this fact throughout the duration of the paper. The next lemma is a stepping stone for Lemma \ref{lem:d1constant}. 

\begin{lemma} \label{lem:MonSquare}
Suppose $p$ is an nc symmetric polynomial that is nc plush on an nc open set, $\cG$. 
Let $2N$ denote the degree of its nc complex hessian, $q$. Then, $q$ must contain 
a term of the form 
$$
\alpha m^T h^T h m \ \ \ (\mbox{or} \ \ \ \alpha m h h^T m^T)
$$ 
where $m$ is an nc analytic monomial of 
degree $N-1$ and $\alpha$ is a positive real constant.
\end{lemma}

\proof Proposition \ref{prop:heredplusantihered} implies $q$ is a sum of 
hereditary and antihereditary polynomials. Let $w$ be a term of degree 
$2N$ in $q$. Without loss of generality, suppose $w$ is hereditary; i.e.,
$w$ has the form 
$$
w = \alpha m_1^T h^T m_2^T m_3 h m_4 
$$
where $\alpha \in \mathbb{R}$, $m_1,m_2,m_3,m_4$ are nc analytic 
monomials in $x$, and 
$$
\deg(m_1) + \deg(m_2) + \deg(m_3) + \deg(m_4) = 2N-2.
$$
By Theorem \ref{thm:P1andP2} (P2), $q$ must contain the Levi-differentially wed term
$$
\tilde{w} = \alpha \tm_1^T h^T h \tm_2
$$
where $\tm_1, \tm_2$ are nc analytic monomials in $x$ and $\deg(\tm_1) = \deg(\tm_2) 
= N-1$.

If $\tm_1 = \tm_2$, we are done (except for showing $\alpha > 0$). 
If the conclusion of the lemma is false, so that $q$ contains no term of the form $\alpha m^T h^T h m$, 
then this implies $\tm_1 \ne \tm_2$. Since $q$ is symmetric, $q$ must also 
contain the term
$$
\tilde{w}^T = \alpha \tm_2^T h^T h \tm_1.
$$
If we partition the border vector so that $e_1^T V = h \tm_1$ and $e_2^T V = h \tm_2$, 
then we get that 
$$
q = \left( \begin{array}{c}
h \tm_1 \\
h \tm_2 \\
\vdots
\end{array} \right)^T
\left( \begin{array}{ccc}
0 & \alpha & \cdots \\
\alpha & 0 & \cdots \\
\vdots & \vdots & \ddots
\end{array} \right) 
\left( \begin{array}{c}
h \tm_1 \\
h \tm_2 \\
\vdots
\end{array} \right).
$$
This middle matrix is not positive semidefinite for any $X \in \cG$. Hence,
Theorem \ref{thm:positivity} implies that $q$ is not positive semidefinite
for all $X \in \cG$. This contradicts the positivity of $q$ on the nc open set, $\cG$. 
Hence, $q$ must contain some term of the form $\alpha m^T h^T h m$. 

We now show $\alpha > 0$. Since we know that $q$ contains a term 
of the form $\alpha m^T h^T h m$ with $m$ an nc analytic or nc antianalytic monomial 
of degree $N-1$, the real constant $\alpha$ will appear on the diagonal in the middle matrix. 
Then, Theorem \ref{thm:positivity} implies that this $\alpha$ must be positive.
\qed 

When we write $e_i$, we mean the vector whose $i^{th}$ entry is $1$ and every other 
entry is $0$. From Equation \eqref{eq:qLDLT}, we can write $q$ as 
a sum of outer products
$$
q = V(x,x^T)[h,h^T]^T \left( \sum_{i=1}^{\cN} {(Le_i)d_i(Le_i)^T} \right) V(x,x^T)[h,h^T] 
$$
\beq  \label{eq:qOuterProd}
= \sum_{i=1}^{\cN}{ V(x,x^T)[h,h^T]^T (Le_i) d_i (Le_i)^T V(x,x^T)[h,h^T]  }.
\eeq
We stress that in Equation \eqref{eq:qOuterProd}, each $Le_i$ and 
$d_i$ depend on $x$ and $x^T$. However, the next lemma shows that one 
element of $D$ is constant and one column of $L$ contains nc polynomials 
rather than nc rationals.

\begin{lemma} \label{lem:d1constant}
Let $p$ be an nc symmetric polynomial that is nc plush on an nc open set. 
Let $2N$ denote the degree of its nc complex hessian, $q$. Then, we can
write the nc complex hessian, $q$, as in Equations \eqref{eq:qLDLT} and 
\eqref{eq:qOuterProd} where $L(x,x^T)$ is unit lower triangular and 
$D(x,x^T) = diag(d_1, \ldots, d_{\cN})$ with $d_1$ a positive real constant.

Hence, each entry in $Le_1$, the first column of $L(x,x^T)$, is an nc polynomial 
rather than an nc rational.
\end{lemma}

\proof Theorem \ref{thm:LDLTCHSY} implies $D(x,x^T)$ is a diagonal matrix. 
Without loss of generality, Lemma \ref{lem:MonSquare} implies that 
$q$ contains a term of the form 
$$
\alpha m^T h^T h m
$$
where $\alpha > 0$ is a positive real constant and $m$ is 
an nc analytic monomial of degree $N-1$. The MMR of $q$ can be written as 
$$
q = \left( \begin{array}{c} 
h m \\
\widehat{V}
\end{array} \right)^T
\left( \begin{array}{cc}
\alpha &  \ell^T \\
\ell & \widehat{M}
\end{array} \right)
\left( \begin{array}{c}
hm \\
\widehat{V}
\end{array} \right).
$$
Since $\alpha > 0$, we can first pivot off $\alpha$ in computing the $LDL^T$ 
factorization of the middle matrix to get
$$
q = 
\left( \begin{array}{c}
hm \\
\widehat{V}
\end{array} \right)^T
\left( \begin{array}{cc}
1 & 0 \\
\frac{1}{\alpha} \ell & I
\end{array} \right)
\left( \begin{array}{cc}
\alpha & 0 \\
0 & \widehat{M} - \frac{1}{\alpha}\ell \ell^T
\end{array} \right)
\left( \begin{array}{cc}
1 & \frac{1}{\alpha}\ell^T \\
0 & I
\end{array} \right)
\left( \begin{array}{c}
hm \\
\widehat{V} 
\end{array} \right).
$$
Now we see that $d_1 = e_1^T D e_1 = \alpha > 0$ and that 
$$
Le_1 = 
\left( \begin{array}{c}
1 \\
\frac{1}{\alpha} \ell
\end{array} \right) \ \ \ \ \ \mbox{and} \ \ \ \ \
\widehat{M} - \frac{1}{\alpha} \ell \ell^T
$$
contain only nc polynomials as entries. \qed 

The next lemma provides even more specific structure to $Le_1$ and 
maintains the nc polynomial structure.

\begin{lemma} \label{lem:LDstructure}
Under the same hypotheses of Lemma \ref{lem:d1constant}, either:
\begin{itemize}
\item[(i)] every entry of $Le_1$ (the 1$^{st}$ column of $L(x,x^T)$) is an nc antianalytic polynomial,
$d_1$ (the 1$^{st}$ diagonal entry of $D(x,x^T)$) is a positive real constant, and the 
corresponding monomials in $V(x,x^T)[h,h^T]$ are nc analytic; or
\item[(ii)] every entry of $Le_1$ (the 1$^{st}$ column of $L(x,x^T)$) is an nc analytic polynomial,
$d_1$ (the 1$^{st}$ diagonal entry of $D(x,x^T)$) is a positive real constant, and the 
corresponding monomials in $V(x,x^T)[h,h^T]$ are nc antianalytic.
\end{itemize}
\end{lemma}

\proof 
Lemma \ref{lem:factorization} implies that $q$ can be written as 
$$
q = A(x)[h]^T Q_1(x,x^T) A(x)[h] + A(x^T)[h^T]^T Q_5(x,x^T) A(x^T)[h^T] 
$$
where each entry of $A(x)[h]$ is an nc analytic monomial and 
each entry of $A(x^T)[h^T]$ is an nc antianalytic monomial. Also, $Q_1$ 
contains hereditary nc polynomials and $Q_5$ contains antihereditary 
nc polynomials. Then, we have that
$$
q = A(x)[h]^T L_1 D_1 L_1^T A(x)[h] + A(x^T)[h^T]^T L_2 D_2 L_2^T A(x^T)[h^T]  
$$
\beq \label{eq:Thm55Eqn}
= 
\underbrace{\left(\begin{array}{c}
A(x)[h] \\
A(x^T)[h^T]
\end{array} \right)^T}_{V(x,x^T)[h,h^T]^T}
\underbrace{\left( \begin{array}{cc}
L_1 & 0 \\
0 & L_2 
\end{array} \right)}_{L}
\underbrace{\left( \begin{array}{cc}
D_1 & 0 \\
0 & D_2 
\end{array} \right)}_{D}
\underbrace{\left( \begin{array}{cc}
L_1 & 0 \\
0 & L_2
\end{array} \right)^T}_{L^T}
\underbrace{\left(\begin{array}{c}
A(x)[h] \\
A(x^T)[h^T]
\end{array} \right)}_{V(x,x^T)[h,h^T]}.
\eeq

Without loss of generality, Lemma \ref{lem:MonSquare} allows us to assume 
that $q$ contains a term of the form 
$$
d_1 m^T h^T h m
$$
where $m$ is an nc analytic monomial in $x$ (so that $hm$ is an entry in 
$A(x)[h]$) of degree $N-1$ and $d_1$ is a positive real constant. Lemma 
\ref{lem:d1constant} implies that $e_1^T D_1 e_1 = d_1$ and that each 
entry of $Le_1$ is an nc polynomial. From Equation \eqref{eq:Thm55Eqn}, 
we have that 
$$
Le_1 = 
\left( \begin{array}{c}
L_1 e_1 \\
0
\end{array} \right)
$$
and $(Le_1)^T V = (L_1 e_1)^T A(x)[h]$. 

Next, write $q$ as in Equation \eqref{eq:qOuterProd} and see that
the first term in this sum becomes
$$
V^T (Le_1) d_1 (Le_1)^T V 
= d_1 ((L_1e_1)^T A(x)[h])^T ((L_1e_1)^T A(x)[h])
$$
Proposition \ref{prop:heredplusantihered} implies that $q$ is a sum of hereditary and 
antihereditary polynomials. Therefore, since $A(x)[h]$ contains only nc analytic 
monomials, this forces $(L_1 e_1)^T$ to contain only nc analytic polynomials 
(which means that $L_1 e_1$ contains only nc antianalytic polynomials). This 
completes the proof of Case (i). 

The proof of Case (ii) works the same way, from Lemma \ref{lem:MonSquare}, 
whenever we assume that $q$ contains a term of the form 
$$
d_1 m h h^T m^T
$$
where $m$ is an nc analytic monomial in $x$ of degree $N-1$ and $d_1$ is a 
positive real constant. \qed 

The next lemma is a technical lemma that is used as a stepping stone to help 
prove Proposition \ref{prop:d1dirder}.

\begin{lemma} \label{lem:d1BV}
Let $p$ be an nc symmetric polynomial that is nc plush on an nc open set. 
Let $2N$ denote the degree of its nc complex hessian, $q$. Then, we can 
write $q$ as in Equation \eqref{eq:qOuterProd}
$$
q = \sum_{i=1}^{\cN}{ V(x,x^T)[h,h^T]^T (Le_i) d_i (Le_i)^T V(x,x^T)[h,h^T] }
$$
with
$$
V(x,x^T)[h,h^T]^T e_1 = x_{i_N}^T \cdots x_{i_2}^T h_{i_1}^T 
\ \ \ \ \ 
(\mbox{resp.} \ \ \ V(x,x^T)[h,h^T]^T e_1 = x_{i_N} \cdots x_{i_2} h_{i_1})
$$
in which case, any term in $q$ that has the form
$$
d_1 \gamma x_{i_N}^T \cdots x_{i_2}^T h_{i_1}^T m(x,h) 
\ \ \ \ \
(\mbox{resp.} \ \ \ d_1 \gamma x_{i_N} \cdots x_{i_2} h_{i_1} m(x^T,h^T) )
$$
where $\gamma$ is a real constant and $m(x,h)$ is some nc analytic monomial 
in $x, h$ of degree 1 in $h$ (resp. $m(x^T, h^T)$ is some nc antianalytic monomial 
in $x^T, h^T$ of degree 1 in $h^T$), is a term in the nc polynomial 
$$
d_1 V(x,x^T)[h,h^T]^T (Le_1) (Le_1)^T V(x,x^T)[h,h^T].
$$ 
Moreover, $\gamma m(x,h)$ (resp. $\gamma m(x^T, h^T)$) is a term in the 
nc analytic (resp. nc antianalytic) polynomial 
$$
(Le_1)^T V(x,x^T)[h,h^T].
$$ 
\end{lemma}

\proof
Proposition \ref{prop:heredplusantihered} implies $q$ is a sum of hereditary and antihereditary 
polynomials. Since the degree of $q$ is $2N$, there exists a term, $w$, in 
$q$ of degree $2N$. Without loss of generality, Lemma \ref{lem:MonSquare} 
allows us to assume that $w$ looks like 
$$
w = d_1 x_{i_N}^T \cdots x_{i_2}^T h_{i_1}^T h_{i_1} x_{i_2} \cdots x_{i_N}.
$$
with $d_1 \in \bbR_+$. We partition the border vector
$V(x,x^T)[h,h^T]$ as
$$
V(x,x^T)[h,h^T] = \left( \begin{array}{c}
h_{i_1} x_{i_2} \cdots x_{i_N} \\
\widehat{V}
\end{array} \right)
$$
where $h_{i_1} x_{i_2} \cdots x_{i_N}$ is not a monomial entry in the vector $\widehat{V}$.
Then, $q$ becomes
\begin{equation*}
q =
\left( \begin{array}{c}
h_{i_1} x_{i_2} \cdots x_{i_N} \\
\widehat{V}
\end{array} \right)^T
\left( \begin{array}{cc}
1 & 0 \\
\ell & \widehat{L}
\end{array} \right)
\left( \begin{array}{cc}
d_1 & 0 \\
0 & \widehat{D}
\end{array} \right)
\left( \begin{array}{cc}
1 & \ell^T \\
0 & \widehat{L}^T
\end{array} \right)
\left( \begin{array}{c}
h_{i_1} x_{i_2} \cdots x_{i_N} \\
\widehat{V}
\end{array} \right)
\end{equation*}
\begin{equation} \label{eq:d1TermsInq}
=
d_1 ( \overbrace{ x_{i_N}^T \cdots x_{i_2}^T h_{i_1}^T h_{i_1} x_{i_2} \cdots x_{i_N} +
x_{i_N}^T \cdots x_{i_2}^T h_{i_1}^T \ell^T \widehat{V} +
\widehat{V}^T \ell h_{i_1} x_{i_2} \cdots x_{i_N}  +
\widehat{V} \ell \ell^T \widehat{V} }^{V(x,x^T)[h,h^T]^T (Le_1)(Le_1)^T V(x,x^T)[h,h^T] =
(x_{i_N}^T \cdots x_{i_2}^T h_{i_1}^T + \widehat{V}^T \ell) (h_{i_1} x_{i_2} \cdots x_{i_N} + \ell^T \widehat{V}) } )
\end{equation}
\begin{equation*}
+ \widehat{V}^T \widehat{L} \widehat{D} \widehat{L}^T \widehat{V}.
\end{equation*}
Since $x_{i_N}^T \cdots x_{i_2}^T h_{i_1}^T$ is not a monomial entry in the vector $\widehat{V}^T$,
this shows that any term in $q$ of the form $d_1 \gamma x_{i_N}^T \cdots x_{i_2}^T h_{i_1}^T m(x,h)$, 
where $\gamma$ is a real constant and $m(x,h)$ is an nc analytic monomial of degree 1 in $h$, 
is a term in the nc polynomial
$$
d_1 V(x,x^T)[h,h^T]^T (Le_1) (Le_1)^T V(x,x^T)[h,h^T].
$$
Equation \eqref{eq:d1TermsInq} implies that either 
$$
d_1 \gamma x_{i_N}^T \cdots x_{i_2}^T h_{i_1}^T m(x,h) = d_1 x_{i_N}^T \cdots x_{i_2}^T h_{i_1}^T h_{i_1} x_{i_2} \cdots x_{i_N} 
$$
or that $d_1 \gamma x_{i_N}^T \cdots x_{i_2}^T h_{i_1}^T m(x,h)$ is a term in the nc polynomial 
$$
d_1 x_{i_N}^T \cdots x_{i_2}^T h_{i_1}^T \ell^T \widehat{V}.
$$
This implies that either $\gamma = 1$ and $m(x,h) = h_{i_1} x_{i_2} \cdots x_{i_N}$ or 
that $\gamma m(x,h)$ is a term in the nc polynomial $\ell^T \widehat{V}$. Hence, 
$\gamma m(x,h)$ is a term in the nc polynomial
$$
h_{i_1} x_{i_2} \cdots x_{i_N} + \ell^T \widehat{V} = (Le_1)^T V(x,x^T)[h,h^T]
$$ 
and Lemma \ref{lem:LDstructure} implies that $(Le_1)^TV(x,x^T)[h,h^T]$ is nc analytic.
\qed

\subsection{1-Differentially Wed Monomials and NC Directional Derivatives} \label{subsec:1dwed}
Earlier, in Section \ref{subsec:Levi}, we introduced Levi-differentially wed monomials and 
provided necessary and sufficient conditions as to when an nc polynomial is an nc complex 
hessian. Now, we introduce \textbf{1-differentially wed monomials} and state Proposition 
\ref{prop:intGvars}, which gives necessary and sufficient conditions as to when a given 
nc polynomial is an nc directional derivative. This machinery is needed to prove Proposition 
\ref{prop:d1dirder}, below.

\subsubsection{Notation} \label{subsec:WedNotation}
Let $m$ be a monomial that is degree one in $h$ or $h^T$. This means $m$ 
contains some $h_i$ or some $h_i^T$. If $m$ contains some $h_i$, we denote 
$$
\monsub{m}{h_i}{x_i}
$$
as the monomial only in $x$ and $x^T$ where $x_i$ replaces $h_i$ in $m$. 
If $m$ contains some $h_i^T$, we denote
$$
\monsub{m}{h_i^T}{x_i^T}
$$
as the monomial only in $x$ and $x^T$ where $x_i^T$ replaces $h_i^T$ in $m$.

Sometimes, we may write $h_i^\gamma$ or $x_i^\gamma$ where $\gamma$ is 
either $\emptyset$ or $T$. When $\gamma = \emptyset$, we define
\begin{eqnarray*}
h_i^\emptyset &:=& h_i \\
x_i^\emptyset &:=& x_i
\end{eqnarray*}
and when $\gamma = T$, we mean $h_i^T$ or $x_i^T$.\\

\begin{example}
If $m = x_1 h_2^T x_1^T$ then $\monsub{m}{h_2^T}{x_2^T} = x_1 x_2^T x_1^T$.
\end{example}

\subsubsection{1-Differentially Wed Monomials} \label{subsec:DiffWedMons}

For $\alpha,\beta$ either $\emptyset$ or $T$, two monomials $m$ 
and $\tm$ are called \textbf{1-differentially wed} if both are degree 
one in $h$ or $h^T$ and if 
$$
\monsub{m}{h_i^\alpha}{x_i^\alpha} = \monsub{\tm}{h_j^\beta}{x_j^\beta}
$$

\begin{example} \label{ex:DiffWed2}
The monomials $m = h_1 x_2^T x_1$ and $\tm = x_1 h_2^T x_1$ are 
1-differentially wed. 
\end{example}

\begin{example} \label{ex:DiffWed3}
The monomials $m = x_2 h_2 x_2$ and $\tm = x_1 x_2 h_2$ are not 
1-differentially wed. 
\end{example}

The following theorem gives necessary and sufficient conditions for an nc 
polynomial to be an nc directional derivative.

\begin{prop} \label{prop:intGvars}
A polynomial $p$ in $x=(x_1,\ldots,x_g)$ and $h=(h_1,\ldots,h_g)$ is 
an nc directional derivative if and only if each monomial in $p$ has degree 
one in $h$ (i.e., contains some $h_j$) and whenever a monomial $m$ 
occurs in $p$, each monomial which is 1-differentially wed to $m$ also 
occurs in $p$ and has the same coefficient.
\end{prop}
\proof This is proved in \cite{GHVpreprint}. \qed 

\subsection{Part of the NC Complex Hessian is an NC Complex Hessian} \label{subsec:PartCHisCH}
In this subsection, we focus on writing the nc complex hessian, $q$, as 
in Equation \eqref{eq:qOuterProd}
$$
q = \sum_{i=1}^{\cN}{ V(x,x^T)[h,h^T]^T (Le_i) d_i (Le_i)^T V(x,x^T)[h,h^T] }.
$$
This subsection culminates with the result that the nc polynomial
$$
d_1 V(x,x^T)[h,h^T]^T(Le_1)(Le_1)^T V(x,x^T)[h,h^T],
$$
is the nc complex hessian for some nc polynomial that is nc plush on 
an nc open set. In order to do this, we first show that the nc polynomial 
$$
(Le_1)^T V(x,x^T)[h,h^T]
$$
is the nc directional derivative for some nc analytic or nc antianalytic polynomial.

\begin{prop} \label{prop:d1dirder}
Let $p$ be an nc symmetric polynomial that is nc plush on an nc open set, $\cG$. 
Let $2N$ denote the degree of its nc complex hessian, $q$. 
If we write $q$ as in Equation \eqref{eq:qOuterProd} and $d_1$ is constant, 
then the nc polynomial
$$
(Le_1)^T V(x,x^T)[h,h^T]
$$
is the nc directional derivative for an nc analytic polynomial or an
nc antianalytic polynomial.

In addition, the nc polynomial
$$
d_1 V(x,x^T)[h,h^T]^T (Le_1)(Le_1)^T V(x,x^T)[h,h^T]
$$
is the nc complex hessian for some nc polynomial that is nc plush on $\cG$.
\end{prop}

\proof Without loss of generality, we can assume, by Lemma \ref{lem:LDstructure}, that
$$
(Le_1)^T V(x,x^T)[h,h^T] 
$$
is an nc analytic polynomial where 
$V(x,x^T)[h,h^T]$ and $Le_1$ are partitioned as 
\begin{equation} \label{eq:BlockLV}
V(x,x^T)[h,h^T] = \left( \begin{array}{c}
V_N \\
V_{N-1} \\
\vdots \\
V_1
\end{array} \right), \ \
Le_1 =
\left( \begin{array}{c}
\ell_0 \\
\ell_1 \\
\vdots \\
\ell_{N-1}
\end{array} \right), \ \
\ell_0 =
\left( \begin{array}{c}
1 \\
\star \\
\vdots \\
\star
\end{array} \right),
\end{equation}
where $\star$ is any nc polynomial and $V_j$ is a vector that contains only nc analytic monomials of the form $h_{i_1} x_{i_2} \cdots x_{i_j}$ having total degree $j$. Each $\ell_j$ is a vector with the same length as $V_j$ and, by Lemma \ref{lem:LDstructure}, $\ell_j$ contains only nc antianalytic polynomials ($\ell_j^T$ contains only nc analytic polynomials). With this setup, we have that
$$
\cF(x,h) := (Le_1)^T V(x,x^T)[h,h^T] = \sum_{j=0}^{N-1} {\ell_j^T V_{N-j}}
$$
is an nc analytic polynomial in $x$ and $h$. We define this as $\cF(x,h)$ for convenience.

Lemma \ref{lem:d1constant} implies $d_1 \in \bbR_+$ is a constant and 
Equation \eqref{eq:qOuterProd} implies that $q$ contains the terms
\begin{equation} \label{eq:qContainsd1FTF}
d_1 V(x,x^T)[h,h^T]^T (Le_1) (Le_1)^T V(x,x^T)[h,h^T] 
= d_1 \left(\sum_{j=0}^{N-1} {\ell_j^T V_{N-j}}\right)^T \left(\sum_{j=0}^{N-1} {\ell_j^T V_{N-j}}\right)
\end{equation}
Then, since the degree of $q$ is $2N$ and the degree of each border vector monomial in $V_{N-j}$ is $N-j$, 
it follows that the degree of each nc analytic polynomial in $\ell_j^T$ is at most $j$.

Lemma \ref{lem:MonSquare} implies that $q$ contains some term of the form
\beq \label{mon:mon1q}
\alpha^2 x_{i_N}^T \cdots x_{i_2}^T h_{i_1}^T h_{i_1} x_{i_2} \cdots x_{i_N}
\eeq
with $\alpha$ a nonzero real constant.
This implies that the vector $V_N$ contains the monomial $h_{i_1} x_{i_2} \cdots x_{i_N}$ 
as an entry. 
Without loss of generality, assume this monomial is first in lexicographic
order. Then, 
$$
e_1^T V(x,x^T)[h,h^T] = e_1^T V_N = h_{i_1} x_{i_2} \cdots x_{i_N}.
$$
As in the proof of Lemma \ref{lem:d1constant}, if $M$ represents the middle matrix 
of $q$, then $e_1^T M e_1 = \alpha^2$ and, after one step in the $LDL^T$ algorithm, 
we see that $\alpha^2 = d_1$. Then, by Theorem \ref{thm:P1andP2} (P2), $q$ also 
contains the Levi-differentially wed terms
\begin{eqnarray*}
&d_1 x_{i_N}^T \cdots x_{i_2}^T h_{i_1}^T x_{i_1} h_{i_2} x_{i_3} \cdots x_{i_N}& \\
&d_1 x_{i_N}^T \cdots x_{i_2}^T h_{i_1}^T x_{i_1} x_{i_2} h_{i_3} \cdots x_{i_N}& \\
& \vdots & \\
& d_1 x_{i_N}^T \cdots x_{i_2}^T h_{i_1}^T x_{i_1} x_{i_2} \cdots x_{i_{N-1}} h_{i_N}.&
\end{eqnarray*}
Since $q$ contains these terms and the term in \eqref{mon:mon1q}, 
Lemma \ref{lem:d1BV} implies that $\cF(x,h)$ contains the term
\begin{equation} \label{mon:MainMonInF}
h_{i_1} x_{i_2} x_{i_3} \cdots x_{i_N}
\end{equation}
and the terms
\begin{eqnarray*}
&x_{i_1} h_{i_2} x_{i_3} \cdots x_{i_N}& \\
&x_{i_1} x_{i_2} h_{i_3} \cdots x_{i_N}& \\
& \vdots & \\
&x_{i_1} x_{i_2} \cdots x_{i_{N-1}} h_{i_N}.&
\end{eqnarray*}
Hence, $\cF(x,h)$ contains all 1-differentially wed monomials to 
$h_{i_1} x_{i_2} x_{i_3} \cdots x_{i_N}$ as terms. Proposition \ref{prop:intGvars} 
implies that $\cF(x,h)$ contains the nc directional derivative of 
$x_{i_1} x_{i_2} x_{i_3} \cdots x_{i_N}$.

Now we pick any other term in $\cF(x,h)$ and show that $\cF(x,h)$ contains
all other 1-differentially wed monomials to it and that they all occur with the 
same coefficient.  Suppose $\cF(x,h)$ contains the term
\begin{equation*}
\gamma x_{s_1} \cdots x_{s_k} h_{\beta_1} x_{\beta_2} \cdots x_{\beta_{N-j}}.
\end{equation*}

We already showed that $\cF(x,h)$ contains the monomial in 
\eqref{mon:MainMonInF}, $h_{i_1} x_{i_2} \cdots x_{i_N}$, 
as a term so $\cF(x,h)^T$ must contain the monomial 
$x_{i_N}^T \cdots x_{i_2}^T h_{i_1}^T$ as a term. This implies that 
$d_1 \cF(x,h)^T \cF(x,h)$ contains the terms
$$
d_1 x_{i_N}^T \cdots x_{i_2}^T h_{i_1}^T 
 (h_{i_1} x_{i_2} \cdots x_{i_N} + 
\gamma x_{s_1} \cdots x_{s_k} h_{\beta_1} x_{\beta_2} \cdots x_{\beta_{N-j}}).
$$
Hence, $q$ contains the term
$$
d_1 \gamma x_{i_N}^T \cdots x_{i_2}^T h_{i_1}^T x_{s_1} \cdots x_{s_k} h_{\beta_1} x_{\beta_2} \cdots x_{\beta_{N-j}}
$$
and Theorem \ref{thm:P1andP2} (P2) implies that $q$ contains the Levi-differentially wed terms
\begin{eqnarray*}
&d_1 \gamma x_{i_N}^T \cdots x_{i_2}^T h_{i_1}^T h_{s_1} x_{s_2} \cdots x_{s_k} x_{\beta_1} x_{\beta_2} \cdots x_{\beta_{N-j}}& \\
&d_1 \gamma x_{i_N}^T \cdots x_{i_2}^T h_{i_1}^T x_{s_1} h_{s_2} \cdots x_{s_k} x_{\beta_1} x_{\beta_2} \cdots x_{\beta_{N-j}}& \\
&\vdots& \\
&d_1 \gamma x_{i_N}^T \cdots x_{i_2}^T h_{i_1}^T x_{s_1} x_{s_2} \cdots h_{s_k} x_{\beta_1} x_{\beta_2} \cdots x_{\beta_{N-j}}& \\
&d_1 \gamma x_{i_N}^T \cdots x_{i_2}^T h_{i_1}^T x_{s_1} x_{s_2} \cdots x_{s_k} x_{\beta_1} h_{\beta_2} \cdots x_{\beta_{N-j}}& \\
&\vdots& \\
&d_1 \gamma x_{i_N}^T \cdots x_{i_2}^T h_{i_1}^T x_{s_1} x_{s_2} \cdots x_{s_k} x_{\beta_1} x_{\beta_2} \cdots h_{\beta_{N-j}}.&
\end{eqnarray*}
Since $q$ contains all of these terms with $x_{i_N}^T \cdots x_{i_2}^T h_{i_1}^T$ on the left,
Lemma \ref{lem:d1BV} implies $\cF(x,h)$ must contain the terms
\begin{eqnarray*}
& \gamma h_{s_1} x_{s_2} \cdots x_{s_k} x_{\beta_1} x_{\beta_2} \cdots x_{\beta_{N-j}}& \\
& \gamma x_{s_1} h_{s_2} \cdots x_{s_k} x_{\beta_1} x_{\beta_2} \cdots x_{\beta_{N-j}}& \\
&\vdots& \\
& \gamma x_{s_1} x_{s_2} \cdots h_{s_k} x_{\beta_1} x_{\beta_2} \cdots x_{\beta_{N-j}}& \\
& \gamma x_{s_1} x_{s_2} \cdots x_{s_k} h_{\beta_1} x_{\beta_2} \cdots x_{\beta_{N-j}}& \\
& \gamma x_{s_1} x_{s_2} \cdots x_{s_k} x_{\beta_1} h_{\beta_2} \cdots x_{\beta_{N-j}}& \\
&\vdots& \\
& \gamma x_{s_1} x_{s_2} \cdots x_{s_k} x_{\beta_1} x_{\beta_2} \cdots h_{\beta_{N-j}}.&
\end{eqnarray*}
All of these terms in $\cF(x,h)$ have the same coefficient, $\gamma$, and are 
1-differentially wed to each other. Thus, Proposition \ref{prop:intGvars} implies 
that they sum to the nc directional derivative of
$$
\gamma x_{s_1} x_{s_2} \cdots x_{s_k} x_{\beta_1} x_{\beta_2} \cdots x_{\beta_{N-j}}.
$$
Hence, we have shown that $\cF(x,h) = (Le_1)^T V(x,x^T)[h,h^T]$ is an nc directional derivative, 
where, without loss of generality, we assumed that $\cF(x,h)$ was nc analytic. 

Now we have that 
$$
\cF(x,h):=(Le_1)^T V(x,x^T)[h,h^T]
$$ 
is the nc directional derivative of some nc analytic or nc antianalytic polynomial. Suppose, 
without loss of generality, that $\cF(x,h)$ is the nc directional derivative of some nc analytic 
polynomial, $\cF(x)$. Then, $\cF(x,h)$ is nc analytic and
$$
d_1 \cF(x,h)^T \cF(x,h) = d_1 V(x,x^T)[h,h^T]^T (Le_1)(Le_1)^T V(x,x^T)[h,h^T]
$$
is the nc complex hessian of the nc polynomial
$$
d_1 \cF(x)^T \cF(x).
$$
Hence, for any $n \ge 1$, any $X \in \cG$, and any $H \in (\bbR^{n \times n})^g$, we have
$$
d_1 \cF(X,H)^T \cF(X,H) \succeq 0.
$$
\qed 

\subsection{Constant $D$ Result} \label{subsec:Dresult}

In this subsection, we show that for an nc symmetric polynomial, $p$, that is
nc plush on an nc open set, the matrix $D(x,x^T)$ in Equation \eqref{eq:qLDLT} 
has no dependence on $x$ or $x^T$ and is actually a positive semidefinite 
constant real matrix. First, we require a helpful lemma.

\begin{lemma} \label{lem:plushDpos}
If $p$ is an nc symmetric polynomial that is nc plush on an nc open set, $\cG$, 
then its nc complex hessian, $q$, can be written as in Equation \eqref{eq:qLDLT}
$$
q = V(x,x^T)[h,h^T]^T L(x,x^T) D(x,x^T) L(x,x^T)^T V(x,x^T)[h,h^T]
$$
where $D(x,x^T)$ is a diagonal matrix of nc rationals and $D(X,X^T) \succeq 0$ for all $X \in \cG$.
\end{lemma}

\proof This follows immediately from Theorem \ref{thm:LDLTCHSY}. \qed

\begin{theorem} \label{thm:ConstantD}
Suppose $p$ is an nc symmetric polynomial that is nc plush on an nc open set, $\cG$. 
Let $2N$ denote the degree of its nc complex hessian, $q$. Then $q$ can be
written as in Equation \eqref{eq:qLDLT}
$$
q = V(x,x^T)[h,h^T]^T L(x,x^T) D(x,x^T) L(x,x^T)^T V(x,x^T)[h,h^T]
$$
where $D(x,x^T) = diag(d_1, d_2, \ldots, d_\cN)$ is a positive semidefinite
constant real matrix (i.e., $d_i \in \bbR_{\ge 0}$ for all $i=1,\ldots,\cN$) and 
$L(x,x^T)$ is a unit lower triangular matrix of nc polynomials.
\end{theorem}

\proof Lemma \ref{lem:plushDpos} implies $D(X,X^T) \succeq 0$ for every
$X \in \cG$. This means $d_i(X,X^T) \succeq 0$
for every $X \in \cG$ and every $i = 1,\ldots, \cN$.
It remains to show that each $d_i$ is a nonnegative constant real number.

First, write the nc complex hessian, $q$, as in Equation \eqref{eq:qOuterProd}
$$
q = \sum_{i=1}^{\cN}{ V(x,x^T)[h,h^T]^T (Le_i) d_i(x,x^T) (Le_i)^T V(x,x^T)[h,h^T]  }.
$$
Lemma \ref{lem:d1constant} shows $d_1 \in \bbR_+$ is a constant, $Le_1$ contains 
nc polynomial entries, and Proposition \ref{prop:d1dirder} shows that
$$
d_1 V(x,x^T)[h,h^T]^T (Le_1) (Le_1)^T V(x,x^T)[h,h^T]
$$
is the nc complex hessian for some nc polynomial that is nc plush on $\cG$. Since nc 
differentiation is linear, we know that the difference of two nc complex hessians 
is an nc complex hessian. This implies that 
\begin{eqnarray*}
\widetilde{q} &:=& q - d_1 V(x,x^T)[h,h^T]^T (Le_1) (Le_1)^T V(x,x^T)[h,h^T] \\
&=& \sum_{i=2}^{\cN}{ V(x,x^T)[h,h^T]^T (Le_i) d_i(x,x^T)(Le_i)^T V(x,x^T)[h,h^T]  }
\end{eqnarray*}
is an nc complex hessian. Since $d_i (X,X^T) \succeq 0$ for all $X \in \cG$
and for all $i$, we have that $\widetilde{q}$ is the nc complex hessian for an nc symmetric
polynomial that is nc plush on $\cG$.
\qed 

Now we give a partial list of references. For a complete list, see \cite{GHVpreprint}.

\end{document}